\documentclass{amsart}

\usepackage{amssymb}
\usepackage{amsmath,amsfonts,amsthm,amstext}

\def\mapnew#1{\smash{\mathop{\longrightarrow}\limits^{#1}}}
\theoremstyle{definition}

\newtheorem{theorem}{Theorem}
\newtheorem{lemma}{Lemma}
\newtheorem{corollary}{Corollary}

\newtheorem{proposition}{Proposition}

\newcommand{\z}{\mathbb{Z}}
\newcommand{\n}{\mathbb{N}}
\newcommand{\real}{\mathbb{R}}

\DeclareMathOperator{\im}{im}

\DeclareMathOperator{\Hom}{Hom}

\DeclareMathOperator{\Aug}{Aug}
\DeclareMathOperator{\Tor}{Tor}

\begin{document}

\author{Laurent Bartholdi}
\thanks{Y. de Cornulier and D. Kochloukova thank the ``Courant
  Research Centre'' of the Georg-August-Universit\"at zu G\"ottingen
  for supporting their visit to G\"ottingen in June 2010 and for
  hospitality during a visit in July 2014.}
\address{
L.B.: Mathematisches Institut\\
Georg-August-Universit\"at zu G\"ottingen\\
D-37073 G\"ottingen, Germany}

\author{Yves de Cornulier}
\thanks{Y. de Cornulier was supported by the
  Brazilian-French Network in Mathematics for the travel funding of
  his visit in Campinas in February, 2010}
\address{
Y.C.: D\'epartement de Math\'ematiques\\
Universit\'e Paris-Sud\\
F-91405 Orsay Cedex, France} 

\author{Dessislava H. Kochloukova}
\thanks{D. Kochloukova is partially supported by ``bolsa de
  produtividade em pesquisa'' from CNPq, Brazil.}
\address{
D.K.: Department of Mathematics\\
State University of Campinas (UNICAMP)\\
CEP 13083-859 Campinas, SP, Brazil}

\date{November 24, 2014}

\title{Homological finiteness properties of wreath products}

\begin{abstract}
  We study the homological finiteness properties $FP_m$ of wreath
  products $\Gamma = H \wr_X G$. We show that, when $H$ has infinite
  abelianization, $\Gamma$ has type $FP_m$ if and only if both $G$ and
  $H$ have type $FP_m$ and $G$ acts (diagonally) on $X^i$ 
with  stabilizers of type $FP_{m-i}$ and with finitely many orbits for all $1
  \leq i \leq m$.

  If furthermore $H$ is torsion-free we give a criterion for $\Gamma$
  to be Bredon-$FP_m$ with respect to the class of finite subgroups of
  $\Gamma$.

  Finally, when $H$ has infinite abelianization and $\chi : \Gamma \to \mathbb{R}$ is a non-zero homomorphism with $\chi(H) = 0$, we classify when $[\chi]$ belongs to
  the Bieri-Neumann-Strebel-Renz invariant $\Sigma^m(\Gamma, \mathbb{Z})$.
\end{abstract}

\maketitle

\section{Introduction}

In this paper we consider the homological finiteness
property $FP_m$ of a group. By definition a group $G$ is of
type $FP_m$ if the trivial $\mathbb{Z} G$-module $\mathbb{Z}$
has a projective resolution with all projectives finitely
generated in dimensions up to $m$. If $G$ is
finitely presented and of type $FP_m$ for some $m \geq 2$
then $G$ is of homotopical type $F_m$ i.e. there is a
$K(G,1)$ with finite $m$-skeleton.

In general it is hard to determine the homological type $FP_m$ of
a group $G$. Even in the case of metabelian groups $G$ there is
no complete classification of the ones of type $FP_m$ for $m > 2$ though
there is an open conjecture, the $FP_m$-Conjecture, that relates the homological type $FP_m$ with the $m$-element subsets of the complement of the Bieri-Strebel invariant $\Sigma_A(Q)$ in the character sphere $S(Q)$ \cite{BieriGroves}. Later on the Bieri-Strebel invariant was generalised for any finitely generated group $G$ \cite{BieriRenz},\cite{BNS} and is widely  referred to as the Bieri-Neumann-Strebel-Renz invariant $\Sigma^1(G, \mathbb{Z})$.  Both properties $FP_2$ and $F_2$ (i.e. finite presentability)  coincide for a metabelian group $G$ and the $FP_2$-Conjecture holds \cite{BieriStrebel}. It is worth mentioning that for a general group $G$ the properties $FP_2$ and $F_2$
do not coincide \cite{BestvinaBrady}. Though there are some sufficient conditions, see \cite{Groves},
there is no 
complete classification of finite presentability or type $FP_2$ even in the
class of  nilpotent-by-abelian groups.

Let $G$, $H$ be groups and $X$ be a $G$-set. The
(permutational) wreath product $\Gamma = H \wr_X G$
is defined to be the semidirect product $H^{(X)} \rtimes G$, where $H^{(X)}$ is the direct sum of $|X|$ copies of $H$ and $G$ acts on $H^{(X)}$  by permuting the summands.
The classification of the
finitely presented wreath products $\Gamma = H \wr_X G$ was
established in \cite{Yves}. It was shown that $\Gamma$ is finitely
presented if and only if both $H$ and $G$ are finitely
presented, $G$ acts  on $X$ with finitely generated stabilizers, and $G$ acts (diagonally) on
$X^2$ with finitely many orbits. We generalize this result by showing in Lemma \ref{m=2} when $G$ is of type $FP_2$. If $H$ has
infinite abelianization we give a
criterion for $\Gamma$ to be of type $FP_m$ for $m \geq 3$. The sufficiency of the conditions of
the criterion do not require that $H$ has infinite
abelianization, see Proposition \ref{easydirection}, but our proof of the necessity of the
conditions uses significantly the fact that $H$ has infinite
abelianization. 

 Our main results are the following theorems. The second is
 a homotopy version of the first one. The proof of Theorem A is homological and Theorem B is an easy corollary of Theorem A and the fact that for $m \geq 2$ 
a group is of type $F_m$ if and only if it is of type $FP_m$ and is finitely presented.

\bigskip
{\bf Theorem A} {\it
Let $\Gamma = H \wr_X G$ be a wreath product, where  $X \not= \emptyset$ and $H$ has infinite abelianization. 
Then the following are equivalent :

1. $\Gamma$  is of type
  $FP_m$;

2. $H$ is of type $FP_m$, $G$ is of type $FP_m$,
$G$ acts on $X^i$ with
  stabilizers of type $FP_{m-i}$ and with finitely many 
orbits   for all $1 \leq i \leq m$.}

\bigskip

{\bf Theorem B}
{\it Let $\Gamma = H \wr_X G$ be a wreath product, where $ X \not= \emptyset$ and $H$ has infinite abelianization. 
Then the following are equivalent :

1. $\Gamma$  is of type
  $F_m$;

2. $H$ is of type $F_m$, $G$ is of type $F_m$,
$G$ acts on $X^i$ with 
  stabilizers of type $FP_{m-i}$ and with finitely many 
orbits  for all $1 \leq i \leq m$.}

\bigskip
As a corollary we obtain some results about the Bredon-$FP_m$ type with respect to the class of finite subgroups, denoted  $\underline{FP}_m$. This requires that $H$ is torsion-free in order to control the conjugacy classes of finite subgroups in $\Gamma = H \wr_X G$. In general a group is of type $\underline{FP}_m$ if it has finitely many conjugacy classes of finite subgroups and the centralizer of every finite subgroup is of type $FP_m$ \cite{KropNucPer}. The homotopical counterpart of the property $\underline{FP}_{\infty}$ was studied earlier in \cite{Luck}.

\bigskip
{\bf Theorem C}
{\it Let  $\Gamma=H\wr_X G$  be a wreath
  product, where $ X \not= \emptyset$ and  $H$ be torsion-free, with
  infinite abelianization.   Then $\Gamma$ has type
  $\underline{FP}_m$ if and only if  the following
  conditions hold :

1. $G$ has type $\underline{FP}_m$;

2. $H$ has type $FP_m$;

3. for every finite subgroup $K$ of $G$ and every $ 1 \leq i
\leq m$ the centralizer
$C_G(K)$ acts on $(K \backslash X )^i$ with
stabilizers of type $FP_{m-i}$ and with finitely many orbits.}

\bigskip
{\it Remark.} If $\Gamma=H\wr_X G$  has type
  $\underline{FP}_m$, $X \not= \emptyset$ and $H$ has a non-trivial finite subgroup $K$ then $X$ is finite. Indeed since in $\Gamma$ there are only finitely many conjugacy classes of finite subgroups there is an upper bound of the order of the finite subgroups. On other hand  for every finite subset $X_i$ of $X$ there is a finite subgroup of $\Gamma$ isomorphic to $K^{X_i}$.
  
  \bigskip
Finally a monoid version of Theorem A is given in Theorem
\ref{sigmamain} for the monoid $\Gamma_{\chi} = \{ g \in
\Gamma | \chi(g) \geq 0 \}$, where $\chi : \Gamma \to
\mathbb{R}$ is a non-trivial character of $\Gamma$
such that $\chi(H) = 0$. This describes the  points $[\chi]$ 
of the Bieri-Neumann-Strebel-Renz invariant $\Sigma^m(\Gamma,
\mathbb{Z})$ with $\chi(H) = 0$ in terms of the invariant $\Sigma^m(G,
\mathbb{Z})$ and the action of $G_{\chi}$ on $X^i$ for $i
\leq m$. 

For a finitely generated group $G$ the Bieri-Neumann-Strebel-Renz
invariants $\{ \Sigma^m(G, \mathbb{Z}) \}_{m \geq 1}$ were first defined in the case $m = 1$  
 in \cite{BNS} and for general $m \geq 1$ were considered in \cite{BieriRenz}.  In general
the homological invariant 
$\Sigma^m(G, \mathbb{Z})$ is an open subset of the unit sphere $S(G)$  and
$\Sigma^m(G, \mathbb{Z})$ determines which subgroups of $G$
above the commutator are of homological type $FP_m$ \cite{BieriRenz}.   
The homological and homotopical $\Sigma^m$-invariants of a group are quite difficult to calculate but they are known for right-angled Artin groups \cite{MMV}, the R. Thompson group $F$ \cite{BGK}, metabelian groups of finite Pr\"ufer rank \cite{Meinert} and for split extensions metabelian groups if $m = 3$ \cite{H-K}.

\subsection*{Structure of the paper} In Section \ref{prel} we collect some basic properties of the homological type $FP_m$ for groups and modules. In section \ref{sectionsemiinduced} we classify the homological type of exterior and tensor powers of induced modules. In Section \ref{split} we classify when a group $G = M \rtimes G$ has type $FP_m$ provided $M$ is a  finitely generated induced $\mathbb{Z} G$-module. The results from Section \ref{sectionsemiinduced} and Section \ref{split} are applied in Section \ref{sectionwreath}, where we prove Theorem A and Theorem B. In Section \ref{sectionbredon} we prove Theorem C and in Section \ref{sigmasection} we prove a $\Sigma$-version of Theorem A.

\section{\boldmath Preliminaries on the homological type $FP_m$}  \label{prel}

Let $R$ be an associative ring with unity.
 We recall that an $R$-module $M$ is of type $FP_m$ if $M$ has a projective resolution with all projectives finitely generated up to dimension $m$. If not otherwise stated $\otimes$ is the tensor product over $\z$ and the modules and the group actions considered are left ones.

\begin{lemma}[{\cite[Prop.~1.4]{BieriQMWbook}}] \label{dimensionshifting1} Let $ 0 \to M_1 \to M  \to M_2 \to 0$ be a short exact sequence of $R$-modules. 

(a) if $M$ and $M_2$ are of type $FP_s$ and $s \geq 1$ then $M_1$ is of type $FP_{s-1}$;

(b) if $M_1$ and $M_2$ are of type $FP_s$ then $M$ is of type $FP_s$;

(c) if $M_1$ is of type $FP_{s-1}$ and $M$ is of type $FP_s$ then $M_2$ is of type $FP_s$.
\end{lemma}

The following results follow easily from Lemma \ref{dimensionshifting1}. For completeness we include proofs.

\begin{lemma}   \label{dimensionshifting2} Let $$ \cdots \to Q_i \mapnew{d_i} Q_{i-1} \to \cdots \to Q_0 \mapnew{d_0} V \to 0$$ be an exact complex of $R$-modules with $Q_i$ of type $FP_{m-i}$ for all $i \leq m$. Then $V$ has type $FP_m$.
\end{lemma}

\begin{proof} Apply Lemma \ref{dimensionshifting1}(c) for the short exact sequences 
$$
0 \to \im (d_{i+1}) = \ker (d_i) \to Q_i \to \im (d_i) \to 0
$$
for $0 \leq i \leq m-1$ to prove by inverse induction on $i$ that $\im (d_i)$ is $FP_{m-i}$.
\end{proof}

\begin{lemma}   \label{dimensionshifting3} Let $ 0 \leq k \leq m$ be integers and let $$ \cdots \to Q_i \mapnew{d_i} Q_{i-1} \to \cdots \to Q_0 \mapnew{d_0} V \to 0$$ be an exact complex of $R$-modules with $Q_i$ of type $FP_{m-i}$ for all $i < k$. Suppose further that $V$ has type $FP_m$. Then $\im (d_k)$ has type $FP_{m-k}$.
\end{lemma}

\begin{proof} Consider the exact complex
$$
0 \to \im (d_k) \to Q_{k-1} \mapnew{d_{k-1}} Q_{k-2} \to \cdots \to Q_0 \mapnew{d_0} V \to 0
$$
and apply Lemma 1(a) for the short exact sequences associated to the above complex.
\end{proof}

\begin{lemma}   \label{dimensionshifting4} Let 
$${\mathcal Q} :  \  \  \cdots \to Q_i \mapnew{d_i} Q_{i-1} \to \cdots \to Q_0 \mapnew{d_0} V \to 0$$ 
be a complex of $R$-modules with $Q_i$ of type $FP_{m-i}$ for all $0 \leq i \leq m$ and $H_i({\mathcal Q})$ of type $FP_{m-i-1}$ for $ 0 \leq i \leq m-1$. Then $V$ has type $FP_m$.
\end{lemma}

\begin{proof} Apply Lemma \ref{dimensionshifting1} for the short exact sequences
$$
0 \to \ker (d_i) \to Q_i \to \im (d_i) \to 0
$$
and
$$
0 \to \im (d_{i+1}) \to \ker (d_i) \to H_i({\mathcal Q}) \to 0
$$ to prove by inverse induction on $i$ that $\im (d_i)$ is $FP_{m-i}$ for $0 \leq i \leq m$.
\end{proof}

The following lemma appeared as a footnote in \cite{Aberg} and was later on explained with more details in \cite[Prop. 4.1]{Bux}.

\begin{lemma} \label{retract} A retract of a group of type $FP_m$ is a group of type $FP_m$ i.e. if the split extension $\Gamma = M \rtimes G$ is of type $FP_m$ then $G$ is of type $FP_m$.
\end{lemma}

Let $\Gamma = M \rtimes G$ be a group. Then $\mathbb{Z} M$ is a left $\mathbb{Z} \Gamma$-module,
where $G$ acts via conjugation on $M$ and $M$ acts on $\z M$ via left multiplication. Note that the augmentation ideal $\Aug(\z M)$ is a $\mathbb{Z} \Gamma$-submodule of $\mathbb{Z} M$.

\begin{lemma} \label{condition} Let $\Gamma = M \rtimes G$
  be a group. Then $\Gamma$ is $FP_m$ if and only if $G$ is $FP_m$ and $\Aug (\mathbb{Z} M)$ is $FP_{m-1}$ as $\mathbb{Z} \Gamma$-module.
\end{lemma}

\begin{proof} By Lemma \ref{dimensionshifting1} applied to the short exact sequence $0 \to \Aug (\z \Gamma) \to \z \Gamma \to \z \to 0$  the group $\Gamma$ is of type $FP_m$ if and only if $\Aug(\z \Gamma)$ is $FP_{m-1}$ as $\z \Gamma$-module.
Consider the short exact sequence of $\z \Gamma$-modules
\begin{equation} \label{ses1}
0 \to \z \Gamma \Aug (\z G) \mapnew{\alpha} \Aug(\z \Gamma) \to \Aug(\z M) \to 0,
\end{equation}
where $\alpha$ is the inclusion map.

If $\Gamma$ is of type $FP_m$ by Lemma \ref{retract} $G$ is of type $FP_m$, so $\Aug(\z G)$ is $FP_{m-1}$ as $\z G$-module and the induced $\z \Gamma$-module  $\z \Gamma \Aug (\z G)$ is $FP_{m-1}$. Then by Lemma 1(c) applied to the short exact sequence
(\ref{ses1}) $\Aug(\z M)$ is $FP_{m-1}$ as $\z \Gamma$-module.

If $G$ is of type $FP_m$ and $\Aug(\z M)$ is $FP_{m-1}$ as $\z \Gamma$-module then by Lemma \ref{dimensionshifting1}(b) applied to the short exact sequence (\ref{ses1}) $\Aug (\z \Gamma)$ is $FP_{m-1}$ as $\z \Gamma$-module.
\end{proof}

\begin{lemma} \label{tensor9} Let $S$ be a subring of $R$ such that $R$ is
  flat as right $S$-module and $M$ be a left $S$-module. Furthermore assume that the
  inclusion map $S \to R$ of right $S$-modules splits.  Then $M$ is of type
  $FP_m$ as $S$-module  if and only if $R \otimes_S M$ is
  $FP_m$ as $R$-module.
\end{lemma}

\begin{proof}
 
1. Suppose that $M$ is of type
  $FP_m$ as $S$-module. 
Let ${\mathcal F}$ be  a projective
  resolution of $M$ with projectives  finitely generated
   in dimensions $\leq m$.
Since $R \otimes_S - $ is an exact
  functor 
$R \otimes_S {\mathcal F}$ is a projective resolution of the $R$-module $R
\otimes_S M$ with projectives finitely generated in
dimensions $\leq m$, so  $R \otimes_S M$ is
  $FP_m$ as $R$-module.

2. Suppose that  $R \otimes_S M$ is
  $FP_m$ as $R$-module. We will prove by induction on $m$ that $M$
  is $FP_m$ as $S$-module.

Suppose first that $m = 0$. Let $X = \{ a_i = \sum_{j}
r_{i,j} \otimes m_{i,j} \}_{i}$
be a finite generating set of  $R \otimes_S M$ as $R$-module,  where $r_{i,j} \in R, m_{i,j} \in M$. Let
$M_0$ be the $S$-submodule of $M$ generated by the finite set $Y = \{
m_{i,j} \}_{i,j}$. Since $R$ is flat as $S$-module
there is a short exact sequence $0 \to R \otimes_S M_0 \mapnew{\alpha}
R \otimes_S M \to  R \otimes_S (M/M_0) \to 0$
and by the choice of $Y$ the map $\alpha$ is surjective, so
$\alpha$ is an
isomorphism. Then $R \otimes_S (M/M_0) = 0$ has a direct
summand $S \otimes_S (M/M_0) = M / M_0$, so $M_0 = M$.

Suppose now that $M$ is $FP_{m-1}$ as $S$-module. Let
$$
{\mathcal F} : \cdots \to F_i \mapnew{d_i} F_{i-1} \to \cdots \to
F_0 \to M \to 0
$$
be a projective resolution of $S$-modules with $F_i$
finitely generated for $i \leq m-1$. Consider the exact
complex induced by ${\mathcal F}$
$$
{\mathcal R} : 0 \to A = \ker (d_{m-1}) \to F_{m-1}
\mapnew{d_{m-1}}  F_{m-2} \to \cdots \to F_0 \to M \to 0.
$$
Since $R$ is flat as right $S$-module 
$$R \otimes_S {\mathcal R}  :    0 \to R \otimes_S A 
\to R \otimes_S F_{m-1}
\mapnew{id_R \otimes d_{m-1}}  R \otimes_S F_{m-2} \to \cdots \to R
\otimes_S F_0 \to R \otimes_S M \to 0
$$ is an
exact complex and $R \otimes_S F_i$ is a  finitely generated
projective $R$-module for $i \leq m-1$. By Lemma \ref{dimensionshifting3} since $R \otimes_S M$ is $FP_m$ as
$R$-module we have that $R \otimes_S A $ is finitely
generated as $R$-module. Then by the case $m = 0$ we deduce
that $A$ is finitely generated as $S$-module, so $M$ is
$FP_m$ as $S$-module.
\end{proof}

\begin{corollary} \label{inductionFPm} Let $L$ be a subgroup
  of a group $G$ and $A$ be a (left) $\z L$-module. Then
 $\z G \otimes_{\z L} A$ is $FP_m$ as $\z G$-module if and only if $A$ is $FP_m$ as $\z L$-module.
\end{corollary}

\section{Exterior and tensor powers of induced modules} \label{sectionsemiinduced}

Let $G$ be a group and $S$ be a free $\mathbb{Z}$-module with basis $X \not= \emptyset$. If $S$ is a $\mathbb{Z} G$-module via an action of  $G$ on the set
$X$, we call $S$ an induced module. If $S$ is a $\mathbb{Z} G$-module via an action of  $G$ on the set
$X \cup - X$ we say that $S$ is a semi-induced module, in this case we write $S = \mathbb{Z}_t X$ with index $t$ standing for twisted action of $G$.

\begin{lemma} \label{twisting} Let $S= \mathbb{Z}_t X$ be a finitely generated semi-induced $\mathbb{Z} G$-module. Then 

(a) $S \simeq \bigoplus_{i \in I} \mathbb{Z} G \otimes_{\mathbb{Z} K_i} V_i$ as $\mathbb{Z} G$-modules where $I$ is finite, for each $i \in I$ there is an element $x_i \in X$ with stabilizer $G_i$ in $G$ such that $K_i$ is a subgroup of $G$ containing $G_i$ such that $[K_i : G_i] = d_i \in \{ 1,2\}$ and each $V_i$ as an abelian group is $\mathbb{Z}$, $G_i$ acts trivially on $V_i$ and $K_i$ acts non-trivially on $V_i$ precisely when $d_i = 2$. 

(b) $S$ is of type $FP_m$ over $\mathbb{Z} G$ if and only if all the groups $G_i$ are of type $FP_m$ for $i \in I$.
\end{lemma}

\begin{proof} (a) Observe first that $S$ splits as a finite direct sum of cyclic semi-induced $\mathbb{Z} G$-modules, each generated by an element of $X$. Thus from the very beginning we may assume that $S$ is a cyclic semi-induced $\mathbb{Z} G$-module with a generator $w$ an element from $X$. If $S$ is an induced $\mathbb{Z} G$-module we are done, so assume that
 $S$  is not an induced $\mathbb{Z} G$-module. 
Then there are $g_1, g_2$ in $G$ such that $g_1 w = - g_2 w$, so for  $t = g_2^{-1} g_1$ we have $t w = - w$. Then $t^2 w = t(-w) = - t w = -(-w) = w$, so $t^2 \in G_0$, where $G_0$ is the stabilizer of $w$ in $G$. For $g_0 \in G_0$ we have $g_0 t w = - g_0 w = - w = t w$, so $t^{-1} g_0 t \in G_0$. Then for the subgroup $K_0$ of $G$ generated by $G_0$ and $t$ we have that $G_0$ is a normal subgroup of $K_0$ of index 2.

Observe that if $h_1, h_2 \in G$ and $h_1 w = - h_2 w$ then for $t_1 = h_2^{-1} h_1$ we have $t_1 w = -w = t w$, so $t_1^{-1} t \in G_0$. Thus $h_1 w = - h_2 w$ is a corollary of $t w = - w$ so $S \simeq (\mathbb{Z} G \otimes_{\mathbb{Z} G_0} \mathbb{Z}) / \mathbb{Z} G (t G_0 + G_0) \simeq \mathbb{Z} G \otimes_{\mathbb{Z} K_0} V_0$, where $V_0 \simeq (\mathbb{Z} K_0 \otimes_{\mathbb{Z} G_0} \mathbb{Z}) / \mathbb{Z} K_0 (G_0 + t G_0)$ has underlying abelian group $\mathbb{Z}$ and $G_0$ acts trivially and $K_0/ G_0$ changes the sign. This completes the proof of item (a).

(b) Observe that since $S \simeq \bigoplus_{i \in I} \mathbb{Z} G \otimes_{\mathbb{Z} K_i} V_i$  we have that $S$ is of type $FP_m$ if and only if
$ \mathbb{Z} G \otimes_{\mathbb{Z} K_i} V_i$ is $FP_m$ as $\mathbb{Z} G$-module for every $i \in I$. By Corollary \ref{inductionFPm} $ \mathbb{Z} G \otimes_{\mathbb{Z} K_i} V_i$ is $FP_m$ as $\mathbb{Z} G$-module precisely when $V_i$ is of type $FP_m$ over $\mathbb{Z} K_i$. If $K_i = G_i$ we have that $V_i$ is the trivial $\mathbb{Z} G_i$-module, so $V_i$ is of type $FP_m$ over $\mathbb{Z} K_i$ precisely when $G_i = K_i$ is $FP_m$ as a group. If $K_i / G_i$ is the cyclic group of order 2 then $V_i$ is the sign module i.e.\ $G_i$ acts trivially on $V_i = \mathbb{Z}$ and $K_i/ G_i$ acts via sign change.  Since $G_i$ has finite index in $K_i$ we have that $V_i$ is $FP_m$ as $\mathbb{Z} K_i$-module precisely when $V_i$ is $FP_m$ as $\mathbb{Z} G_i$-module, but as $\mathbb{Z} G_i$-module $V_i$ is the trivial one, so $V_i$ is $FP_m$ as $\mathbb{Z} G_i$-module is the same as $G_i$ is $FP_m$. 
\end{proof}

 The following proposition is the main result of this section. It classifies the homological type of tensor and exterior powers of induced modules.

\begin{proposition} \label{wedgeproducts} Let $M= \mathbb{Z} X$ be an induced $\mathbb{Z} G$-module. Then the following are equivalent :

1.  $\bigwedge^i M$ is of type $FP_{m-i}$ for all $1 \leq i \leq m$;

2.   $\bigotimes^i M$ is of type $FP_{m-i}$ for all $1 \leq i \leq m$;

3. all stabilizers of the diagonal action of $G$ on $X^i$   are of type $FP_{m-i}$  and
 $\otimes^i M$ is finitely generated as  $\mathbb{Z} G$-module for all $1 \leq i \leq m$;

4. all stabilizers  of $i$ element subsets of $X$  are of type $FP_{m-i}$  and
 $\wedge^i M$ is finitely generated as  $\mathbb{Z} G$-module for all $1 \leq i \leq m$;
 
5.  all stabilizers of the diagonal action of $G$ on $X^i$  are of type $FP_{m-i}$  and $G \backslash X^i $ is finite for all $1 \leq i \leq m$.

\end{proposition}

\begin{proof}
Observe that $\wedge^i M = \z_t Y_i$ is an semi-induced $\mathbb{Z} G$-module, where $Y_i$ is the set of $i$-element subsets of $X$ and $\otimes^i M = \z Z_i$ is an induced $\mathbb{Z} G$-module, where $Z_i = X^i$.

It is obvious that conditions 2,3 and 5  are equivalent.

Note that the stabilizer of $(x_1, \ldots, x_i) \in X^i$ in $G$ (via the diagonal action)  is  commensurable with the stabilizer of the set $\{ x_1, \ldots, x_i \}$ in $G$, in particular one of the stabilizers is $FP_{m-i}$ precisely when the other is $FP_{m-i}$.
 Furthermore $\wedge^i M$ is finitely generated as $\mathbb{Z}
 G$-module for all $i \leq m$ if and only if $G$ acts on the set of $i$ element subsets
 of $X$ with finitely many orbits for all $i \leq m$. The last is equivalent with  $G$ acts on $X^i$ with finitely many orbits for all $i \leq m$.
Thus conditions 5 and 4 are equivalent, and by Lemma \ref{twisting}(b) conditions 1 and 2 are equivalent too.
\end{proof}  

\begin{lemma} \label{wedgeproducts3} 
Let $X$ be a (left) $G$-set such that  $G$ acts  with finitely many orbits on $X^i$ and 
with stabilizers  of type $FP_{m-i}$ for all $1 \leq i \leq m$. 
Let $i_1, i_2, \ldots, i_s$ be non-negative integers such that $i_1 + \cdots + i_s = r \leq m$ 
 and
$X^{(i_1,i_2, \ldots,i_s)}
$ be the set of $s$-tuples $(X_1,
\ldots, X_s)$ of
pairwise disjoint subsets $X_1, \ldots , X_s$ of $X$ such
that $X_j$ has cardinality $i_j$ for all $1 \leq j \leq s$.
Then $\mathbb{Z} X^{(i_1,i_2, \ldots,i_s)}$ is of type $FP_{m-r}$ as $\mathbb{Z} G$-module. 
\end{lemma}

\begin{proof}
Let $B$ be the stabilizer in $G$ of the $s$-tuple $(X_1,
\ldots, X_s)$ i.e. $B$ is the set of $g \in G$ such that $g X_i = X_i $ for all $ 1 \leq i \leq s$. Then since all $X_i$ are finite $B$ is a subgroup of finite index in the group $D$, where $D$ is  the stabilizer  in $G$ of the set $\bigcup_{1 \leq i \leq s} X_i$. Note that $| \bigcup_i X_i | = i_1 + \cdots + i_s = r$. Then by Proposition \ref{wedgeproducts} (i.e. item 5 implies item 4), $D$ is of type $FP_{m-r}$, hence $B$ is of type $FP_{m-r}$.
Note that $G \backslash X^r $ finite implies that $G \backslash X^{(i_1,i_2, \ldots,i_s)}$ is finite. Then by Lemma \ref{twisting} the  $\mathbb{Z} G$-module $\mathbb{Z} X^{(i_1,i_2, \ldots,i_s)}$ is of type $FP_{m-r}$.
\end{proof}

\section{\boldmath Type $FP_m$ for split extensions of induced modules} \label{split}

Let $\Gamma = M \rtimes G$ be a group, where $M$ is abelian.
Note we do
not impose restrictions on $G$. We consider $M$ as a left $\z G$-module, where $G$ acts by conjugation. In this section $M$ is an induced finitely generated (left) $\mathbb{Z} G$-module i.e. \begin{equation} \label{directsum} M = \bigoplus_{i \in I} \mathbb{Z} G  \otimes_{\mathbb{Z} H_i}\mathbb{Z} \simeq \bigoplus_{i \in I} \mathbb{Z} [G/ H_i] = \mathbb{Z} X,
 \end{equation}
where $X = \bigcup_{i \in I} G/ H_i$, 
 $I$ is a finite set and $\{ H_i \}_{i \in I}$ are the stabilizers of the action of $G$ on $X$. In this section we classify when a group $\Gamma = M \rtimes G$ is of type $FP_m$ provided $M$  is an induced finitely generated $\mathbb{Z} G$-module.

\begin{proposition} \label{free1}  Let $G$ be a group of
  type $FP_m$ and $M$ be a finitely generated induced $\mathbb{Z}[G]$-module such that  $\wedge^i M$ be of type $FP_{m-i}$ as $\mathbb{Z} G$-module via the diagonal $G$-action for all $1 \leq i \leq m$.  Then $\Gamma = M \rtimes G$ is of type $FP_m$.
\end{proposition}

\begin{proof}
Let $X$ be the disjoint union $\bigcup_i G/ H_i$, so $X$ is a
basis of $M$ as a free $\mathbb{Z}$-module.
Since $M$ is a torsion-free abelian group, there is a Koszul complex \cite{Weibel}
\begin{equation} \label{abelian1}
\cdots \to \mathbb{Z} M \otimes_{\mathbb{Z}} (\wedge^k M) \mapnew{d_k}  \mathbb{Z} M \otimes_{\mathbb{Z}} (\wedge^{k-1} M) \to \cdots \to \mathbb{Z} M \otimes_{\mathbb{Z}} M \to \mathbb{Z} M \to \mathbb{Z} \to 0
\end{equation}
with differential given by 
$$
d_k(m_1 \wedge \cdots \wedge m_k) = \sum_{1 \leq i \leq k} (-1)^i \epsilon(m_i) \otimes m_1 \wedge \cdots \wedge m_{i-1} \wedge m_{i+1} \wedge \cdots \wedge m_k,
$$
where $\epsilon(m_i) = m_i - 1 \in \z M$ and $m_1, \ldots,
m_k \in X$.
Thus (\ref{abelian1}) gives an exact complex
\begin{equation} \label{abelian1new}
{\mathcal Q} : \cdots \to Q_{k-1} = \mathbb{Z} M \otimes_{\mathbb{Z}} (\wedge^k M) \mapnew{d_k}  Q_{k-2} = \mathbb{Z} M \otimes_{\mathbb{Z}} (\wedge^{k-1} M) \to \cdots $$ $$ \to Q_0 = \mathbb{Z} M \otimes_{\mathbb{Z}} M \to Q_{-1} = \Aug(\mathbb{Z} M) \to 0.
\end{equation}
Note that for a left $\mathbb{Z} G$-module $V$ there is an isomorphism of abelian groups
$$
\mathbb{Z} \Gamma \otimes_{\mathbb{Z} G } V \simeq \mathbb{Z} M \otimes_{\mathbb{Z}} V
$$
and thus the action of $\Gamma$ on the induced module $\mathbb{Z} \Gamma \otimes_{\mathbb{Z} G } V $ gives an action of $\Gamma$ on $\mathbb{Z} M \otimes_{\mathbb{Z}} V$.
We apply this for $V = \wedge^k M$ and since $\wedge^k M$ is of type $FP_{m-k}$ as $\mathbb{Z} G$-module and  
$\mathbb{Z} \Gamma \otimes_{\mathbb{Z} G }$ is an exact functor  
$$
\mathbb{Z} \Gamma \otimes_{\mathbb{Z} G } (\wedge^k  M) \simeq \mathbb{Z} M \otimes_{\mathbb{Z}} (\wedge^k M)
$$
is of type $FP_{m-k}$ as $\mathbb{Z} \Gamma$-module. 
Note that since $X$ is a $G$-invariant set the differentials in (\ref{abelian1new}) are homomorphisms of $\z G$-modules, hence are homomorphisms of $\z \Gamma$-modules.
Then (\ref{abelian1new}) is an exact complex with $Q_i$ of homological type $FP_{m-i-1}$ as $\mathbb{Z} \Gamma$-module for $0 \leq  i \leq m-1$. Then by Lemma \ref{dimensionshifting2} applied to the complex  (\ref{abelian1new}) $\Aug (\mathbb{Z} M)$ is of type $FP_{m-1}$ as $\mathbb{Z} \Gamma$-module, hence by Lemma \ref{condition} $\Gamma$ is of type $FP_m$.
\end{proof}

\begin{lemma} \label{free2} Suppose that $\Gamma = M \rtimes
  G$ is a group of type $FP_m$ and $M$ is a finitely generated induced  $\mathbb{Z}[G]$-module such that $\wedge^i M$ is of type $FP_{m-i-1}$ as $\mathbb{Z} G$-module (via the diagonal $G$-action) for all $1 \leq  i \leq m -1$. Then  $M$ is of type $FP_{m-1}$ as $\mathbb{Z} G$-module.
\end{lemma}

\begin{proof}  By Lemma \ref{condition} $\Omega = \Aug(\mathbb{Z} M)$ is of type $FP_{m-1}$ as $\mathbb{Z} \Gamma$-module. Note that $M \simeq \Omega / \Omega^2 \simeq \mathbb{Z} \otimes_{\mathbb{Z} M} \Aug (\mathbb{Z} M)$. Let
$$
{\mathcal F} : \cdots \to F_i \to F_{i-1} \to \cdots \to F_0 \to \Omega \to 0
$$
be a free resolution of $\Omega$ as $\mathbb{Z} \Gamma$-module with $F_i$ finitely generated for $i \leq m-1$. 
Let ${\mathcal Q}$ be the complex $\mathbb{Z} \otimes_{\mathbb{Z} M} {\mathcal F}$ i.e.
 \begin{equation} \label{complexQ}
{\mathcal Q} : \cdots \to Q_i \to Q_{i-1} \to \cdots \to Q_0 \to M \to 0
\end{equation}
is a complex of free $\mathbb{Z} G$-modules with $Q_i$ finitely generated for $i \leq m-1$ and its homology groups
are $$
H_i({\mathcal Q}) \simeq \Tor_{i}^{\mathbb{Z} M}(\mathbb{Z}, \Omega) \hbox{ for } i \geq 1.
$$
Note that the long exact sequence in $\Tor$ applied to the short exact sequence $\Omega \to \mathbb{Z} M \to \mathbb{Z}$ gives 
$$
\Tor_{i}^{\mathbb{Z} M}(\mathbb{Z}, \Omega) \simeq \Tor_{i+1}^{\mathbb{Z} M}(\mathbb{Z}, \mathbb{Z}) = H_{i+1}(M, \mathbb{Z}) \simeq \wedge^{i+1} M \hbox{ for } i \geq 1.
$$
The last isomorphism is \cite[Ch.~V,~Thm.~6.4]{Brownbook}.
Then $H_i({\mathcal Q})$ is of type $FP_{m-i-2}$ as $\mathbb{Z} G$-module for all $ i \leq m-2$. Then by Lemma \ref{dimensionshifting4} applied to the complex (\ref{complexQ})
$M$ is of type $FP_{m-1}$ as $\mathbb{Z} G$-module.
\end{proof}

\begin{proposition} \label{free3} Let $\Gamma = M \rtimes G$
  be a group of type $FP_m$ and  
$M$ be a finitely generated induced $\mathbb{Z}[G]$-module. Then $\wedge^i M$ is of type $FP_{m-i}$ as $\mathbb{Z} G$-module via the diagonal $G$-action for all $1 \leq i \leq m$.
\end{proposition}

\begin{proof} We induct on $m$ and assume that the proposition holds for smaller values of $m$, in particular $\wedge^i M$ is $FP_{m-i-1}$ as $\mathbb{Z} G$-module. By Lemma \ref{free2} $M$ is of type $FP_{m-1}$ as $\mathbb{Z} G$-module, thus the proposition holds for $i = 1$. We proceed by induction on $i$ i.e. assume that we have proved  that $\wedge^j M$ is of type $FP_{m-j}$ as $\mathbb{Z} G$-module via the diagonal $G$-action for all $1 \leq j \leq i-1$ and will  show that 
$\wedge^i M$ is of type $FP_{m-i}$ as $\mathbb{Z} G$-module. 

Consider the Koszul complex 
\begin{equation} \label{abelian2}
\cdots \to \mathbb{Z} M \otimes_{\mathbb{Z}} (\wedge^k M) \mapnew{d_k}  \mathbb{Z} M \otimes_{\mathbb{Z}} (\wedge^{k-1} M) \to \cdots \to \mathbb{Z} M \otimes_{\mathbb{Z}} M \mapnew{d_1} \mathbb{Z} M \to \mathbb{Z} \to 0
\end{equation}
and its modified version
\begin{equation} \label{abelian*new}
{\mathcal Q} : \cdots \to Q_{k-1} = \mathbb{Z} M \otimes_{\mathbb{Z}} (\wedge^k M) \mapnew{d_k}  Q_{k-2} = \mathbb{Z} M \otimes_{\mathbb{Z}} (\wedge^{k-1} M) \to \cdots $$ $$\to Q_0 = \mathbb{Z} M \otimes_{\mathbb{Z}} M \mapnew{d_1} Q_{-1} = \Aug(\mathbb{Z} M) \to 0
\end{equation}
Since $\wedge^k M$ is $FP_{m-k}$ as $\mathbb{Z} G$-module for $k \leq i-1$ the induced module
$Q_{k-1} = \mathbb{Z} M \otimes_{\mathbb{Z}} (\wedge^k M)$ is $FP_{m-k}$ as $\mathbb{Z} \Gamma$-module for $k \leq i-1$. Note that by Lemma \ref{condition} $\Aug (\z M)$ is $FP_{m-1}$ as $\z \Gamma$-module.
Then by Lemma \ref{dimensionshifting3} applied to the complex (\ref{abelian*new}) $\im (d_i)$ is of type $FP_{m-1-(i-1)}$ as $\mathbb{Z} \Gamma$-module.
Note that
$$\im (d_{i}) \simeq (\mathbb{Z} M \otimes (\wedge^i M)) / \ker (d_{i}) \simeq
(\mathbb{Z} M \otimes (\wedge^i M)) / \im (d_{i+1}).$$
Denote $(\mathbb{Z} M \otimes (\wedge^i M)) / \im (d_{i+1})$ by $V$ i.e. $V$ has type $FP_{m-i}$ as $\z \Gamma$-module. 
Let $W$ be $\mathbb{Z} \otimes_{\mathbb{Z} M} V$. Since $ \mathbb{Z} \otimes_{\mathbb{Z} M} -$ is right exact there is
an exact sequence
$$
\mathbb{Z} \otimes_{\mathbb{Z} M} \im (d_{i+1}) \mapnew{\beta} \mathbb{Z} \otimes_{\mathbb{Z} M} (\mathbb{Z} M \otimes (\wedge^{i} M)) \to W \to 0.$$
By the definition of $d_{i+1}$ we have $\beta = 0$, furthermore the module in the middle is isomorphic to $\wedge^i M$, so
$$
W \simeq \wedge^i M.$$ If $W$ is of type $FP_{m-i}$ as $\mathbb{Z} G$-module the inductive step is completed.

Let
$$
{\mathcal F} : \cdots \to F_j \to F_{j-1} \to \cdots \to F_0 \to V \to 0
$$
be a free resolution of $V$ as $\mathbb{Z} \Gamma$-module with $F_j$ finitely generated for $j \leq m-i$.  
  Let ${\mathcal R}$ be the complex $\mathbb{Z} \otimes_{\mathbb{Z} M} {\mathcal F}$ i.e.
 $$
{\mathcal R} : \cdots \to R_j \to R_{j-1} \to \cdots \to R_0 \to W \to 0
$$
is a complex of free $\mathbb{Z} G$-modules with $R_j$ finitely generated for $j \leq m-i$ and its homology groups are
$$
H_j({\mathcal R}) \simeq \Tor_{j}^{\mathbb{Z} M}(\mathbb{Z}, V) \hbox{ for } j \geq 1.
$$
Note that for the $\z M$-module  $V = \im (d_i)$ the Koszul complex (\ref{abelian2}) gives a free resolution 
\begin{equation} 
{\mathcal C} : \cdots \to \mathbb{Z} M \otimes_{\mathbb{Z}} (\wedge^k M) \mapnew{d_k}  \mathbb{Z} M \otimes_{\mathbb{Z}} (\wedge^{k-1} M) \to \cdots \to $$ $$\mathbb{Z} M \otimes_{\mathbb{Z}} (\wedge^{i+1} M) \to  \mathbb{Z} M \otimes_{\mathbb{Z}} (\wedge^{i} M) \mapnew{d_i} V \to 0
\end{equation}
and by the definition of the differentials $d_k$ the complex 
\begin{equation}  \mathbb{Z} \otimes_{\mathbb{Z} M} {\mathcal C} :
 \cdots \to  \wedge^{i+2} M \to  \wedge^{i+1} M \to  \wedge^{i} M \mapnew{\alpha} \mathbb{Z} \otimes_{\mathbb{Z} M} \im (d_i) = W \to 0
\end{equation}
has all zero differentials except possibly $\alpha$.
Thus 
$$
H_j({\mathcal R}) \simeq \Tor_{j}^{\mathbb{Z} M}(\mathbb{Z}, V) \simeq H_j(\mathbb{Z} \otimes_{\mathbb{Z} M} {\mathcal C}) \simeq \wedge^{j+i} M  \hbox{ for } j \geq 1.
$$
By Lemma \ref{dimensionshifting4} applied to the complex ${\mathcal R}$, for $W$ to be of type $FP_{m-i}$ as $\z G$-module it is sufficient that
$$H_j({\mathcal R}) \hbox{ is of type } FP_{m-i-j-1} \hbox{ as } \mathbb{Z} G\hbox{-module for } 1 \leq j \leq m-i-1.$$ Thus we need for $k = i+j$ that 
$$\wedge^k M \hbox{ is of type } FP_{m-k-1} \hbox{ as } \mathbb{Z} G\hbox{-module for } 1 \leq k \leq m-1$$
but this follows from the fact that the proposition holds for $m-1$.
\end{proof}

The following theorem is the main result of this section and classifies when  a group $\Gamma = M \rtimes G$ is of type $FP_m$ provided $M$  is an induced finitely generated $\mathbb{Z} G$-module.

\begin{theorem} \label{mainthm1}  Let $\Gamma = M \rtimes G$
  be a group and $M$ be a finitely generated induced $\mathbb{Z}[G]$-module. Then
  $\Gamma$ is of type $FP_{m}$ if and only if $G$ is of type
  $FP_m$ and $\wedge^i M$ is of type $FP_{m-i}$ as $\z
  G$-module via the diagonal action for all $1 \leq i \leq m$.
\end{theorem}

\begin{proof} Follows directly from Lemma \ref{retract}, Proposition \ref{free1} and Proposition \ref{free3}.
\end{proof}
 
\section{\boldmath Homological type $FP_m$ for wreath products} \label{sectionwreath}

\subsection{Preliminaries on tensor products of complexes} \label{sectionKunneth}
Let $({\mathcal A}, \partial)$ be a non-negative complex of free
$\z$-modules i.e. all non-zero modules are in
dimension $\geq 0$. Let ${\mathcal A}_i = {\mathcal A}$ for all $1 \leq i \leq s$. 
Consider the tensor product ${\mathcal A}_1 \otimes
\cdots \otimes {\mathcal A}_s$ with
differential  given by
\begin{equation} \label{differential01}
d(a_1 \otimes \cdots \otimes a_s) = 
\sum_{1 \leq j \leq s}
(-1)^{\deg(a_1) + \cdots + \deg(a_j)} a_1 \otimes \cdots
\otimes \partial(a_j)\otimes \cdots \otimes a_s.
\end{equation}
For the transposition $\pi = (i, i+1) \in S_s$ and each
$a_i$ an element of one of the free modules in ${\mathcal A}_i$ we define
\begin{equation} \label{symmetricgroup}
\pi(a_1 \otimes \cdots \otimes a_s) := (-1)^{\deg(a_i) \deg(a_{i+1})} (a_1 \otimes \cdots \otimes  a_{i+1} \otimes a_i \otimes \cdots \otimes a_s).
\end{equation}
Since the symmetric group $S_s$ is generated by $\{ (i, i+1) \}_{1 \leq i \leq
  s-1}$ we get an action of $S_s$ on  
the tensor product ${\mathcal A}_1 \otimes
\cdots \otimes {\mathcal A}_s$. As the action of the symmetric
group $S_s$ commutes with the differential from (\ref{differential01}), this induces an
action of $S_s$ on the homology groups of ${\mathcal A}_1 \otimes
\cdots \otimes {\mathcal A}_s$.

\subsection{Wreath products : proofs of Theorem A and Theorem B}

In this section we classify when a wreath product  $\Gamma = H \wr_X G$  is of type $FP_m$. First we show  in Proposition 
\ref{easydirection} some sufficient conditions for $\Gamma$ to be of type $FP_m$. The difficult part of the classification is to show that these sufficient conditions are necessary. We establish this in Theorem \ref{2direction} under the extra condition that $H$ has infinite abelianization.

\begin{proposition} \label{easydirection} Let $\Gamma = H \wr_X G$ be a wreath
  product such that both $H$ and $G$ are of type $FP_m$,
  $G$ acts (diagonally) on $X^i$  with 
  stabilizers of type $FP_{m-i}$ and with finitely many orbits for $1 \leq i \leq m$. Then
  $\Gamma$ has type $FP_m$.
\end{proposition}

\begin{proof}
Let 
$$
{\mathcal H} :  \cdots \to \mathbb{Z} H \otimes \mathbb{Z} Y_i \mapnew{\partial_i}
\mathbb{Z} H \otimes \mathbb{Z} Y_{i-1} \to \cdots \to
\mathbb{Z} H \otimes \mathbb{Z} Y_1 \mapnew{\partial_1} \mathbb{Z} H \to
\mathbb{Z} \to 0
$$
be a free resolution of the trivial $\mathbb{Z} H$-module
$\mathbb{Z}$ with all $Y_i$ finite  for $i \leq
m$.

Consider $H_x$ an isomorphic copy of the group $H$ and let
${\mathcal H}_x$  be the  complex obtained from ${\mathcal H}$ by
substituting $H$ with $H_x$. We write $\partial_{x}$ for the
differential of ${\mathcal H}_x$ and if we want to stress its
degree, say $i$, we write $\partial_{i,x}$. Let $M$ be the normal closure of $H$ in $\Gamma$. Thus $M$
is the subgroup of the direct product $\prod_{x \in X} H_x$ that
contains the elements with all but finitely
many trivial coordinates. 

We can take the tensor product over $\mathbb{Z}$ of the deleted complexes ${\mathcal H}_x^{del}$ (i.e. in ${\mathcal H}_x$ delete $\mathbb{Z}$) for $x \in X$ 
and obtain the deleted complex ${\mathcal F}^{del}$.  
This makes sense for infinite $X$ as the direct limit of the
tensor products of any finite number of the complexes ${\mathcal
  H}_x$. Note that we have fixed a linear order $\leq$ on $X$
and the tensor product ${\mathcal H}^{del}_{x_1} \otimes \cdots
\otimes {\mathcal H}^{del}_{x_t}$ is made for $x_1 < \cdots <
x_t$ and $t \geq 1$.
By the K\"unneth formula
 we get a free resolution of the trivial $\mathbb{Z}
 M$-module $\mathbb{Z}$
$$
{\mathcal F} : \cdots \to \mathbb{Z} M \otimes W_i \mapnew{d_i}
\mathbb{Z} M \otimes W_{i-1} \to \cdots
\to \mathbb{Z} M \otimes W_1 \mapnew{d_1} \mathbb{Z} M \to
\mathbb{Z} \to 0,
$$
where $W_i$ is a free abelian group with a basis $Z_i$, $Z_i$ is the disjoint union
$$
Z_i = \bigcup (Y_1^{i_1}) \times
(Y_2^{i_2} ) \times \cdots \times (Y_s^{i_s}) \times X^{(i_1,i_2, \ldots,i_s)},
$$
over all possible $s \geq 1$ such that  $i_1 +2i_2
+ \cdots+s i_s = i$, $i_j \geq 0$ and
$
X^{(i_1,i_2, \ldots,i_s)}$ is the set of $s$-tuples $(X_1,
\ldots, X_s)$ of
pairwise disjoint subsets $X_1, \ldots , X_s$ of $X$ such
that $X_j$ has cardinality $i_j$ for all $1 \leq j \leq s$.
We write the elements of $Z_i$ as formal products $y_1
\cdots y_j x_1 \cdots x_j$, which indicates that in the tensor product $y_i \in Y
= \bigcup_{i \geq 1} Y_i$ is taken from the complex ${\mathcal
  H}_{x_i}$ and $x_1 < \cdots < x_j$ are elements of $X$,
recall we have fixed some linear order $\leq$ on $X$. Indeed the element $y_1
\cdots y_j x_1 \cdots x_j$ corresponds to the element  $B_1 \times \cdots \times B_s \times (X_1, \ldots , X_s)$ of $Z_i$, where $X_r$ is the set of those elements $x_t$ of $\{ x_1, \ldots , x_j \}$ for which $y_t \in Y_r$, we write $X_r = \{ x_{j_1} , \ldots , x_{j_{i_r}} \}$ where $x_{j_1} < \cdots < x_{j_{i_r}}$ and set $B_r = (y_{j_1}, \ldots , y_{j_{i_r}}) \in Y_r^{i_r}$.  

Sometimes it will be convenient to use general products $y_1
\cdots y_j x_1 \cdots x_j$, where $x_1, \ldots , x_j$ are pairwise different elements of $X$ but we do not assume that $x_1 < \cdots < x_j$.  We follow (\ref{symmetricgroup}) and define for all $ 1 \leq i \leq j-1$
$$ y_1
\cdots y_j x_1 \cdots x_j = (-1)^{\deg(y_i) \deg (y_{i+1})} y_1
\cdots y_{i+1} y_i \cdots y_j x_1 \cdots x_{i+1} x_i \cdots x_j.
$$
This completes the definition of the general products $y_1
\cdots y_j x_1 \cdots x_j$ and shows that they belong to $Z_i \cup - Z_i$.

As $X$ is a (left) $G$-set we get an action of $G$ on the general products given by
\begin{equation} \label{actionaction}
g (y_1
\cdots y_j x_1 \cdots x_j) = y_1
\cdots y_j (g x_1) \cdots (g x_j).\end{equation}
Thus $W_i$ is a semi-induced $\mathbb{Z} G$-module i.e. $W_i = \mathbb{Z}_t Z_i$. Furthermore since $Y_i$ is a finite set for $i \leq m$ and by Lemma \ref{wedgeproducts3} we deduce that $W_i$ is finitely generated as $\mathbb{Z} G$-module.
Using the other notation for the elements of $Z_i$ we have for
$B_1 \times \cdots \times B_s \times (X_1, \ldots , X_s) \in Z_i$ and $g \in G$ that
$$
g (B_1 \times \cdots \times B_s \times (X_1, \ldots , X_s)) = \pm 
(g B_1 \times \cdots \times g B_s \times (g X_1, \ldots , g X_s))
$$
where $g B_k$ is obtained from $B_k \in Y_k^{i_k}$ by some permutation (possible the trivial one) of the $i_k$ coordinates. By Lemma 
  \ref{twisting}(b) $W_i$ is of type $FP_{m-i}$ as $\mathbb{Z} G$-module if the stabilizer in $G$ of the element $B_1 \times \cdots \times B_s \times (X_1, \ldots , X_s) \in Z_i$ is of type $FP_{m-i}$ for all  $B_1 \times \cdots \times B_s \times (X_1, \ldots , X_s) \in Z_i$ . Any such stabilizer is commensurable with the stabilizer in $G$ of  $ (X_1, \ldots , X_s) \in X^{(i_1,i_2, \ldots,i_s)}$ and by the proof of Lemma \ref{wedgeproducts3}  the last stabilizer is of type $FP_{m-i}$. 
Note we have proved that $W_i$ is of type $FP_{m-i}$ as $\mathbb{Z} G$-module for $1
\leq i \leq m$. 
Then $\mathbb{Z} \Gamma \otimes_{\mathbb{Z} G} W_i \simeq  \mathbb{Z} M \otimes W_i$ is of type $FP_{m-i}$ as $\mathbb{Z} \Gamma$-module. 

Since the action of the symmetric group in (\ref{symmetricgroup}) commutes with the differential 
(\ref{differential01}) for a general product $y_1 \cdots y_j x_1 \cdots x_j$
the differential of
 ${\mathcal F}$ is 
\begin{equation} \label{maindif}
d_s(y_1 \cdots y_j x_1 \cdots x_j)
= (\sum_{\deg(y_i) \geq 2} (-1)^{\deg(y_1) + \cdots + \deg(y_i)}   y_1 \cdots \partial_{x_i}(y_i) \cdots y_j x_1
\cdots x_i \cdots x_j) +$$ $$ (\sum_{\deg(y_i) =1}   (-1)^{\deg(y_1) + \cdots + \deg(y_i)}    y_1 \cdots \partial_{x_i}(y_i)
\cdots y_j x_1 \cdots \widehat{x_i} \cdots x_j)
\end{equation}
where $s = \sum_{1  \leq i  \leq j}  \deg(y_i) \geq 1$. Note that by (\ref{actionaction}) the above differential commutes with the $G$-action, hence commutes with the $\Gamma$-action.
Then by Lemma \ref{dimensionshifting2} 
applied to the complex (obtained from ${\mathcal F}$)
\begin{equation} \label{+}
\widetilde{\mathcal F} : \cdots \to \mathbb{Z} M \otimes W_i \mapnew{d_i}
\mathbb{Z} M \otimes W_{i-1} \to \cdots
\to \mathbb{Z} M \otimes W_1 \mapnew{d_1} \Aug( \mathbb{Z} M) \to 0,
\end{equation}
we deduce that
$\Aug(\mathbb{Z} M)$ is of type $FP_{m-1}$ as $\mathbb{Z} \Gamma$-module, so by Lemma  \ref{condition}
$\Gamma$ is of type $FP_m$.
\end{proof}
 
\begin{lemma} \label{m=2} Let $\Gamma = H \wr_X G$ be a group with $H \not= 1, X \not= \emptyset$. Then
  $\Gamma$ is of type $FP_2$ if and only if  $H$ is of type
  $FP_2$, $G$ is of type $FP_2$ and 
 $G$ acts on $X^i$  with
  stabilizers of type $FP_{2-i}$ and with finitely many orbits  for $1 \leq i \leq 2$.
\end{lemma}

\begin{proof} 
By the previous
theorem it remains to show that if $\Gamma$ is $FP_2$ then all stated conditions hold. This is a particular case of \cite[Prop.~4.9]{BCS} when $X$ is $G$-transitive, the general case follows by  the same argument.
\end{proof}

\begin{proposition} \label{dificultdirection} Let $m \geq
  1$, $\Gamma = H \wr_X G$ be a group of type
  $FP_m$ such that $H \not= 1, X \not= \emptyset$, $H$ is of type $FP_{m-1}$ and 
 $G$ acts on $X^i$  with
  stabilizers of type $FP_{m-i}$ and with finitely many orbits  for $1 \leq i \leq m$.
Then $H$ is of type $FP_m$.
\end{proposition}

\begin{proof}
As pointed in \cite{Yves} the case $m = 1$ holds, The case $m = 2$ is a particular case of Lemma \ref{m=2}. So
we may assume that $m \geq 3$.
 
 Let $Z_i$ and $Y_i$ be as in the proof of
  Proposition \ref{easydirection}. 
Since $H$ is of type $FP_{m-1}$   
  we  can assume that
  $$Y_1,Y_2, \ldots, Y_{m-1} \hbox{ are all finite }$$ and
will prove that $Y_m$ can be chosen finite.
 The description of $Z_{m}$ and $X^{(0, \ldots, 0,1)} = X$ imply
 that
$$\mathbb{Z} Z_{m} = T_{1,m} \oplus T_{2,m},$$ where $T_{1,m} =
\mathbb{Z} [Y_m \times X]$ and $T_{2,m} = \mathbb{Z} [Z_m
\setminus (Y_m \times X)]$ i.e. $T_{1,m} =
  \mathbb{Z} Y_m \otimes \mathbb{Z} X$. 
Observe that since $Y_1, \ldots, Y_{m-1}$ are all finite and 
 $G$ acts on $X^i$ with finitely many orbits  for $1 \leq i \leq m$ we have that $T_{2,m}$ is finitely generated as
$\mathbb{Z} G$-module. 
Note that by the proof of Proposition \ref{easydirection} $\mathbb{Z} M \otimes \mathbb{Z} Z_i$ is $FP_{m-i}$ as $\mathbb{Z} \Gamma$-module for $1 \leq i \leq m-1$.
Since $\Gamma$ is of type $FP_m$ by Lemma \ref{condition} and 
Lemma \ref{dimensionshifting3} applied to the complex (\ref{+})   $\im (d_m)$ is finitely
generated as $\mathbb{Z} \Gamma$-module. Consider the
splitting of the domain $\mathbb{Z} M \otimes \mathbb{Z}
Z_m$ of the differential $d_m$ as $$\mathbb{Z} M \otimes \mathbb{Z}
Z_m = W_{1,m} \oplus W_{2,m},$$ where $W_{1,m} = \mathbb{Z} M \otimes T_{1,m}$
and $W_{2,m} = \mathbb{Z} M \otimes T_{2,m}$. Since
$W_{2,m}$ is finitely generated as $\mathbb{Z}
\Gamma$-module the
condition that $\im (d_m)$ is finitely generated is
equivalent with $d_m(W_{1,m}) / d_m(W_{1,m}) \cap
d_m(W_{2,m})$ is finitely generated as $\mathbb{Z} \Gamma$-module.

Let $w_1 \in
W_{1,m}, w_2 \in W_{2,m}$ be such that $d_m(w_1) =
d_m(w_2)$. Then $w_1 - w_2 \in \ker (d_m) = \im (d_{m+1})$.
Let $$p : W_{1,m} \oplus W_{2,m} \to W_{1,m}$$ be the
canonical projection.
By the definition of $W_{1,m}$ and (\ref{maindif}) applied with $s =
m +1$ we get that
$$w_1 = p(w_1 - w_2) \in p( \im (d_{m+1})) = $$ $$  p d_{m+1}((\mathbb{Z} M \otimes
\mathbb{Z} [Y_1 \times  Y_m \times X^{(1,0, \ldots, 0, 1)}])
\oplus (\mathbb{Z} M \otimes
\mathbb{Z} [Y_{m+1} \times X])
).$$
Then since $d_{m+1}(\mathbb{Z} M \otimes
\mathbb{Z} [Y_{m+1} \times X]
) \subseteq W_{1,m}$ we have $p d_{m+1}(\mathbb{Z} M \otimes
\mathbb{Z} [Y_{m+1} \times X]) = d_{m+1}(\mathbb{Z} M \otimes
\mathbb{Z} [Y_{m+1} \times X])$ and so
$$d_m(w_1) \in  d_m p d_{m+1}(\mathbb{Z} M \otimes
\mathbb{Z} [Y_1 \times  Y_m \times X^{(1,0, \ldots, 0, 1)}])
+ d_m p d_{m+1}(\mathbb{Z} M \otimes
\mathbb{Z} [Y_{m+1} \times X]
) = $$ $$d_m p d_{m+1}(\mathbb{Z} M \otimes
\mathbb{Z} [Y_1 \times  Y_m \times X^{(1,0, \ldots, 0, 1)}])
+ d_m  d_{m+1}(\mathbb{Z} M \otimes
\mathbb{Z} [Y_{m+1} \times X]) = $$ 
\begin{equation} \label{eqnovo5}  d_m p d_{m+1}(\mathbb{Z} M \otimes
\mathbb{Z} [Y_1 \times  Y_m \times X^{(1,0, \ldots, 0, 1)}]).\end{equation}
Recall that $\partial_{j,x}$ denotes the differential of
${\mathcal H}_x$ in dimension $j$ and when we do not want to
stress the dimension we omit the index $j$ and  use $\partial_x$.
Note that for $y_1 \in Y_1, y_m \in Y_m, x_1, x_m \in X, x_1 \not= x_m$
$$
d_m p d_{m+1} (y_1 y_m x_1 x_m) = d_m p(- \partial_{x_1} (y_1) y_m
x_m  + (-1)^{m+1} y_1 \partial_{x_m}(y_m) x_1 x_m) = 
$$
\begin{equation} \label{eqwr0} d_m(-  \partial_{x_1} (y_1) y_m
x_m) = (-1)^{m+1} \partial_{x_1} (y_1) \partial_{x_m}(y_m) x_m. \end{equation}
Decompose $M = S_x \times H_x$, where $S_x
\subseteq  \prod_{t \not= x} H_{t}$. Thus in (\ref{eqwr0}) we have   
$\partial_{x_m}(y_m) \in \im (\partial_{m,x_m}),  \partial_{x_1} (y_1) \in \Aug (\mathbb{Z} S_{x_m})$.
Then by (\ref{eqnovo5}) and (\ref{eqwr0}) 
$$
d_m(w_1) \in \bigoplus_{y_m \in Y_m,y_1 \in Y_1, x_1, x_m \in
  X, x_1 \not= x_m} \mathbb{Z}
M\partial_{x_1}(y_1) \partial_{x_m}(y_m) x_m
=$$
\begin{equation} \label{eqwr1} \bigoplus_{x_m \in X} \Aug(\mathbb{Z} S_{x_m}) \mathbb{Z}
H_{x_m} \im (\partial_{m,x_m}) x_m = 
\bigoplus_{x\in X} \Aug(\mathbb{Z}
S_{x}) 
\im (\partial_{m,x}) x =: J
\end{equation}
and
$$ d_m(W_{1,m}) = \bigoplus_{x \in X} \mathbb{Z} M
 \im (\partial_{m,x}) x =$$ 
\begin{equation} \label{eqwr2} 
\bigoplus_{x \in X} (\mathbb{Z} S_x  \mathbb{Z}
H_x) \im (\partial_{m,x}) x = \bigoplus_{x \in X} \mathbb{Z} S_x  \im (\partial_{m,x}) x.
\end{equation}
Furthermore for every $j \in J$ there is $w_1 \in W_{1,m}$ such that there is some $w_2 \in W_{2, m}$ with $d_m(w_1) = d_m(w_2)$ and $d_m(w_1) = j$.
Thus by (\ref{eqwr1}) and (\ref{eqwr2})
$$K := d_m(W_{1,m}) /J \simeq  \bigoplus_{x \in X}  \im
(\partial_{m,x}) x 
$$
Thus the condition that $\im (d_m)$
is finitely generated as $\mathbb{Z} \Gamma$-module is
equivalent with 
$K$ is finitely generated as $\mathbb{Z}
\Gamma$-module. Take a finite set $D$ of generators of
$K$ such that $D
\subseteq \bigcup_{x \in X} \im (\partial_{m,x}) x$. Write $d
\in D$ as
$ a_x x$ and define $A_x \subseteq \im (\partial_{m,x})$ as the set of all possible
$a_x$. Thus the set of the images of the finite set $\bigcup_{x \in X}
A_x$ in $\im (\partial_m)$, where we send canonically every
$H_x$ to $H$ by just forgetting the index $x$,  generates $\im (\partial_m)$ as $\mathbb{Z} H$-module. Thus $\im
(\partial_m)$ is finitely generated as $\mathbb{Z}
H$-module, so $Y_m$ can be chosen finite.

\end{proof}

\begin{theorem} \label{2direction} Let $m \geq 1$ and $\Gamma = H \wr_X G$ be of type
  $FP_m$, where $ X \not= \emptyset$ and $H$ has infinite abelianization. Then $H$ is of type $FP_m$
and  $G$ acts on $X^i$ with 
  stabilizers of type $FP_{m-i}$ and with finitely many 
orbits  for all $1 \leq i \leq m$.
\end{theorem}

\begin{proof} We induct on $m$. The case $m = 1$ is
  trivial as the property $FP_1$ is equivalent with finite
  generation and then  \cite[Prop.~2.1]{Yves} applies.

Since $H$ has infinite abelianization $\mathbb{Z}$ is a
retract of $H$, so $\Gamma_0 =\mathbb{Z} \wr_X G$ is a retract of
$\Gamma$. By Lemma \ref{retract} since $\Gamma$ is of type $FP_m$ the group
$\Gamma_0$ is $FP_m$ too. Then by Proposition \ref{wedgeproducts} and Theorem \ref{mainthm1} $G$ acts on $X^i$ with finitely many 
orbits  and 
  stabilizers of type $FP_{m-i}$ for all $1 \leq i \leq m$.
 Then Proposition \ref{dificultdirection} completes the
inductive step.
\end{proof}

\begin{theorem} \label{mainwreathresult}
Let $\Gamma = H \wr_X G$ be a wreath product, where $ X \not= \emptyset$ and $H$ has infinite abelianization. 
Then the following are equivalent :

1. $\Gamma$  is of type
  $FP_m$;

2. $H$ is of type $FP_m$, $G$ is of type $FP_m$,
$G$ acts on $X^i$ with 
  stabilizers of type $FP_{m-i}$ and with finitely many 
orbits  for all $1 \leq i \leq m$.
\end{theorem}
 
\begin{proof}
Follows from Lemma \ref{retract}, Proposition
\ref{easydirection} and Theorem \ref{2direction}.
\end{proof}

The following is a homotopy version of Theorem \ref{mainwreathresult}.

\begin{theorem} \label{mainwreathresult2}
Let $\Gamma = H \wr_X G$ be a wreath product, where $X \not= \emptyset$ and $H$ has infinite abelianization. 
Then the following are equivalent :

1. $\Gamma$  is of type
  $F_m$;

2. $H$ is of type $F_m$, $G$ is of type $F_m$,
$G$ acts on $X^i$ with 
  stabilizers of type $FP_{m-i}$ and with finitely many 
orbits  for all $1 \leq i \leq m$.
\end{theorem}
 Remark. Observe that in the condition 2 we have a mixture
 of both properties $F_m$ and $FP_m$ and not just the
 homotopical property $F_m$.

\begin{proof} In the case $m \leq 2$ this is the main result of
  \cite{Yves}. So we can assume that $m > 2$. In general a
  group is of type $F_m$ if and only if it is finitely
  presented (i.e. is $F_2$) and is of type $FP_m$. Then the
  result follows from Theorem \ref{mainwreathresult} and the fact that the case $m = 2$ of Theorem \ref{mainwreathresult2} was already proven  in  \cite{Yves}.
\end{proof}

 \section{An example}
 
Let $F$ be the Richard Thompson group with infinite presentation $$\langle x_0, x_1, x_2 \ldots | x_j^{x_i} =  x_{j+1}  \hbox{ for } 0 \leq i < j \rangle.$$ Consider its realization as piecewise linear transformations of the interval $[0,1]$, see \cite{essay}. Let $X =(0, 1) \cap \z [\frac{1} {2}]$. Then $F$ acts with finitely many orbits on $X^i$ for every $i \geq 1$ i.e. the $F$ orbit of $x = (x_1, \ldots, x_i) \in X^i$ contains $y = (y_1, \ldots , y_i)$ if and only if there is a permutation $\sigma \in S_i$ such that $x_{\sigma(1)} \leq x_{\sigma(2)} \leq \cdots \leq x_{\sigma(i)}$ and   
 $y_{\sigma(1)} \leq y_{\sigma(2)} \leq \cdots \leq y_{\sigma(i)}$ and $x_{\sigma(k)} = x_{\sigma(k+1)}$ if and only if $y_{\sigma(k)} = y_{\sigma(k+1)}$. Furthermore the stabilizer in $F$ of the point  $x = (x_1, \ldots, x_i) \in X^i$ is $F^{j+1}$, where $j$ is the number of different coordinates of $x$. By \cite{GeogheganBrown} $F$ is of type $FP_{\infty}$, so all stabilizers of the action of $F$ on $X^i$ are of type $FP_{\infty}$. The following result is a direct corollary of Theorem \ref{mainwreathresult} and Theorem \ref{mainwreathresult2}.
 
 \begin{corollary} \label{wreathF} Let $X$ and $F$ be as
   above  and $\Gamma = H \wr_X F$ be a wreath product. Then
   $H$ has type $FP_m$ (resp. $F_m$) if and only if $\Gamma$ is of type $FP_{m}$ (resp. $F_m$).
 \end{corollary} 
 
\begin{proof}
Follows directly from Proposition \ref{easydirection} and Proposition \ref{dificultdirection}.
\end{proof}
 
\section{\boldmath Bredon type $FP_m$ for wreath products} \label{sectionbredon}

\subsection{Preliminaries on Bredon homology}

 In this section we study the homological finiteness Bredon type  $\underline{FP}_m$. Its homotopical $\infty$-version, the homotopical finiteness Bredon property $\underline{F}_{\infty}$ was considered
in \cite{Luck}. Here we need the following homological
version of L\"uck's result from \cite{Luck}.

\begin{theorem} \label{Bredonref}  \cite{KropNucPer}  A group $\Gamma$ is of type
  $\underline{FP}_m$ if and only if $\Gamma$ has finitely
  many conjugacy classes of finite subgroups and for every
  finite subgroup $K$ of $\Gamma$ the
  centralizer $C_{\Gamma}(K)$ is of type $FP_m$.
\end{theorem}

\subsection{\boldmath Bredon type $FP_m$ for wreath products  : a
  proof of Theorem C }

\begin{theorem} \label{finite}  Let $\Gamma=H\wr_X G$ be a wreath product such that $H$ is torsion-free.
Then every finite subgroup of $\Gamma$ is conjugate to a subgroup of
$G$.
\end{theorem}

\begin{proof} In this proof we consider $X$ as a right $G$-set (if
 $X$ is equipped with a left $G$-action it becomes
a right $G$-action by $x g = g^{-1} x$). 
Denote by $\pi : \Gamma\to G$ the natural homomorphism. Let $M$ be the normal closure of $H$ in $\Gamma$. Elements
  $\gamma$ of $\Gamma$ may uniquely be written as $\gamma=mg$ for
  $m\in M$ and $g\in G$; we have $\pi(\gamma)=g$, and write
  $\gamma@x:=m_x$, where $m_x$ is the $x$-coordinate of $m
  \in M \subseteq H^X$.
  The assertion of the theorem is that for every
  finite subgroup $K\le\Gamma$ there exists $\delta\in\Gamma$ such
  that $(k^\delta)@x=1$ for all $x\in X$, $k \in K$.

  Consider first two points $x,y\in X$ in the same orbit under
  $\pi(K)$. Then there exists a unique $h_{x,y}\in H$ such that
  $k@x=h_{x,y}$ for all $k\in K$ with $x^{\pi(k)}=y$. Indeed, consider
  two such $k,k'\in K$ i.e. $x^{\pi(k)}=y = x^{\pi(k')}$. Then $\pi(k'k^{-1})$ fixes $x$, and has finite
  order $s$ because $k'k^{-1}$ belongs to $K$; so
  $1 = (k'k^{-1})^s@x=((k'k^{-1})@x)^s$ implies
  $1 = (k'k^{-1})@x=(k'@x)(k@x)^{-1}$. Note the cocycle identity
  $h_{x,y}h_{y,z}=h_{x,z}$.

  In each orbit $\Omega = x^{\pi(K)}$ of $\pi(K)$ on $X$, choose a representative
  $x_\Omega$, and define $\delta\in H^X$ by $\delta_x=h_{x,x_\Omega}$
  whenever $x$ lies in the orbit $\Omega$.
  Note that since $K$ is finite  there are only finitely
  many $x,y\in X$ such that $h_{x,y}\neq1$. Therefore $\delta$ is a
  finitely supported function on $X$
  i.e. $\delta \in M$.

  Now consider $k\in K$, and write $y=x^{\pi(k)}$. Note that
  $x_{\Omega} = y_{\Omega}$. Then 
  \[(k^\delta)@x=(\delta^{-1}@x)(k@x)(\delta@y)=h_{x,x_\Omega}^{-1}h_{x,y}h_{y,y_\Omega}=h_{x,x_\Omega}^{-1}h_{x,y_\Omega}
  = 1.\qedhere\]
\end{proof}


\begin{theorem} Let  $\Gamma=H\wr_X G$  be a wreath
  product, $ X \not= \emptyset$  and $H$ be torsion-free, with
  infinite abelianization.   Then $\Gamma$ has type
  $\underline{FP}_m$ if and only if  the following
  conditions hold :

1. $G$ has type $\underline{FP}_m$;

2. $H$ has type $FP_m$;

3. for every finite subgroup $K$ of $G$ and every $ 1 \leq i
\leq m$ the centralizer
$C_G(K)$ acts on $(K \backslash X)^i$ with stabilizers of type $FP_{m-i}$ and with  finitely many orbits.
\end{theorem}

\begin{proof} By Theorem \ref{Bredonref} we have to understand the centralizers and the conjugacy classes of finite subgroups in $\Gamma$. By Theorem \ref{finite}
  every finite subgroup $K$ of $\Gamma$ is conjugated to a finite
  subgroup $K_0$ of $G$. Thus $\Gamma$ has finitely many
  conjugacy classes of finite subgroups if and only if $G$
  has finitely many conjugacy classes of finite
  subgroups. This completes the proof for $m = 0$.

Furthermore for a finite subgroup $K$ of $G$ we have that $$C_{\Gamma}(K) = C_M(K) \rtimes C_G(K) \simeq H
\wr_{K \backslash X} C_G(K),$$
where $M$ is the normal closure of $H$ in $\Gamma$ and we consider $X$ as a left $G$-set. Then by Theorem \ref{mainwreathresult}
the centralizer $C_{\Gamma}(K) $ is of type $FP_m$ if and only if $C_G(K)$ is $FP_m$ and conditions 2 and  3 hold. 
Then  Theorem \ref{Bredonref} completes the proof.
\end{proof}

 \section{\boldmath $\Sigma$-invariants} \label{sigmasection}

\subsection{Some basic properties}
 Following
\cite{BieriRenz} we define an equivalence relation $\sim$ on $\Hom(G, \mathbb{R}) \setminus \{ 0 \}$ by $\chi_1 \sim \chi_2$ if there is $r \in \mathbb{R}_{>0}$ such that $\chi_1 = r \chi_2$. Denote by $[\chi]$ the equivalence class of $\chi$. Then the character sphere $S(G)$ is $(\Hom(G, \mathbb{R}) \setminus \{ 0 \}) / \sim$. By definition
$$
\Sigma^m(G, \mathbb{Z} ) = \{ [\chi] \in S(G)  | \  \ \mathbb{Z} \hbox{ is }
FP_m \hbox{ as } \mathbb{Z} G_{\chi}-\hbox{module} \},
$$where
$G_{\chi} = \{ g \in G | \chi(g) \geq 0 \}$.

\begin{lemma} \label{sigmaretract}  Let $m \geq 1$ be an integer, $\Gamma = M \rtimes G$ be a finitely generated group  and $[\chi] \in \Sigma^m(\Gamma, \z)$ such that $\chi(M) = 0$. Then
$[\chi_0] \in \Sigma^m (G, \z)$, where $\chi_0$ is the restriction of $\chi$ to $G$.
\end{lemma}

\begin{proof} This is a monoid version of Lemma \ref{retract} and Aberg's idea  ( i.e. the footnote from \cite{Aberg}) can be easily modified in this context.
Indeed by \cite{BieriQMWbook} we have that $[\chi] \in \Sigma^m(\Gamma, \z)$ is equivalent with 
$$
\Tor^{\z \Gamma_{\chi}}_i(\z, \prod \z \Gamma_{\chi}) = 0 \hbox{ for all }1 \leq i \leq m-1
$$
and the trivial $\mathbb{Z} \Gamma_{\chi}$-module $\mathbb{Z}$ is finitely presented. Note that by the functoriality of $\Tor_i^{\z S} ( \z , \prod \z S)$ on $S$ applied to the monoids $S \in \{ G_{\chi_0}, \Gamma_{\chi} \}$ and the fact that $\Gamma_{\chi} = M \rtimes G_{\chi_0}$ we get that  $\Tor^{\z G_{\chi_0}}_i(\z, \prod \z G_{\chi_0}) = 0$ for all $1 \leq i \leq m-1$. Furthermore the trivial $\mathbb{Z} G_{\chi_0}$-module $\mathbb{Z}$
is obtained from the trivial $\mathbb{Z} \Gamma_{\chi}$-module $\mathbb{Z}$  applying the  right exact functor $\mathbb{Z} \otimes_{\mathbb{Z} M} $, hence the trivial $\mathbb{Z} G_{\chi_0}$-module $\mathbb{Z}$ is finitely presented.
Thus $[\chi_0] \in \Sigma^m (G, \z)$.
\end{proof}

\begin{lemma} \label{condition2} Let $m \geq 1$ be an integer, $\Gamma =  M \rtimes G$ 
be a finitely generated group
and $\chi : \Gamma \to \real$ be a non-zero homomorphism such that $\chi(M) = 0$. 
Denote by $\chi_0$ the restriction of $\chi$ to $G$.
Then $[\chi] \in \Sigma^m(\Gamma, \z)$ if and only if $[\chi_0] \in \Sigma^m( G, \z)$ and $\Aug (\mathbb{Z} M)$ is $FP_{m-1}$ as $\mathbb{Z} \Gamma_{\chi}$-module.
\end{lemma}

\begin{proof} In the proof of Lemma \ref{condition} substitute $\Gamma$ with $\Gamma_{\chi}$ and substitute $G$ with $G_{\chi_0}$.
\end{proof}

 The next lemma collects results from \cite{Master}, that can be found in English version in \cite{MMV}.

\begin{lemma} \label{sigmaofactions} \cite{Master}, \cite{MMV} Let $k \geq 0$ be an integer, $\chi : G \to \real$ be a non-zero  homomorphism, $L$ a subgroup of $G$. Denote by $\chi_0$ the restriction of $\chi$ on $L$. Then 

(1) $\z G_{\chi}$ is a flat (right or left) $\z L_{\chi_0}$-module and the inclusion map of $\z L_{\chi_0}$-modules $\z L_{\chi_0} \to \z G_{\chi}$ splits;

(2) if $\chi_0 \not= 0$ then for every $\mathbb{Z} L$-module $V$ we have an isomorphism of  left $\mathbb{Z} G_{\chi}$-modules $\mathbb{Z} G \otimes_{\mathbb{Z} L} V \simeq \mathbb{Z} G_{\chi} \otimes_{\mathbb{Z} L_{\chi_0}} V$;

(3) if $\chi_0 \not= 0$ then for every $\mathbb{Z} L$-module $V$ we have that $\mathbb{Z} G \otimes_{\mathbb{Z} L} V$ is $FP_m$ as $\mathbb{Z} G_{\chi}$-module if and only if $V$ is $FP_m$ as $\mathbb{Z} L_{\chi_0}$-module;

(4) if $L$ has finite index in $G$ then for every $\mathbb{Z} G$-module $V$ we have that $V$ is $FP_m$ as $\mathbb{Z} G_{\chi}$-module if and only if
$V$ is $FP_m$ as $\mathbb{Z} L_{\chi_0}$-module.
\end{lemma}

\begin{proof}
Consider
 $T$  a left transversal of $L$ in $G$  such that $1 \in T$. Then
$
\z G_{\chi} = \bigoplus_{t \in T} t \z L_{\chi_0 \geq - \chi(t)},
$
where $ L_{\chi_0 \geq r} = \{ g \in L | \chi_0(g) \geq r \}$.
Since $1 \in T$ the inclusion map of (right) $\z L_{\chi_0}$-modules $\z L_{\chi_0} \to \z G_{\chi}$ splits.

The rest of item 1. is \cite[Lemma~9.1(i)]{MMV}. Item 2 is \cite[Lemma~9.1(ii)]{MMV}, item 3 is \cite[Lemma~9.2]{MMV} and item 4 is \cite[Thm~9.3]{MMV}.
\end{proof}

 \subsection{Semi-induced modules}
 
 \begin{lemma} \label{twisting2} Let $\chi: G \to \real$ be
   a non-zero homomorphism. 
Let $S$ be a finitely generated semi-induced $\mathbb{Z}
G$-module and consider the decomposition  $S \simeq \bigoplus_{i \in I} \mathbb{Z} G \otimes_{\mathbb{Z} K_i} V_i$ given by Lemma \ref{twisting}. Thus there is a subgroup $G_i$ of $K_i$ of index at most 2 that acts trivially on $V_i$. Then

1. $S$ is finitely generated as $\mathbb{Z} G_{\chi}$-module  if and only if for all $i \in I$ for the restriction $\chi_i$ of $\chi$ on $K_i$ we have $\chi_i \not= 0$;

2. $S$ is of type $FP_m$ as $\mathbb{Z} G_{\chi}$-module if and only if for every $i \in I$ for the restriction $\tilde{\chi}_i$ of $\chi$ to $G_i$ we have that $\widetilde{\chi}_i \not= 0$ and the trivial $\mathbb{Z} (G_i)_{\tilde{\chi}_i}$-module $\mathbb{Z}$ has type $FP_m$ i.e.\ $[\tilde{\chi}_i] \in \Sigma^ m(G_i, \mathbb{Z})$.
\end{lemma}

\begin{proof} Note that item 1 is obvious. By item 1 we can assume that $\widetilde{\chi}_i \not= 0$. 
By Lemma \ref{sigmaofactions} (3) $\mathbb{Z} G \otimes_{\mathbb{Z} K_i} V_i $ is $FP_m$ as $\mathbb{Z} G_{\chi}$-module if and only if $V_i$ is $FP_m$ as $\mathbb{Z} (K_i)_{\chi_i}$-module. 

It remains to consider only the case when $K_i \not= G_i$, so $K_i / G_i$ is cyclic of order 2.
Observe that $V_i$ is the sign module with $G_i$ acting trivially. By Lemma \ref{sigmaofactions} (4) $V_i$ is $FP_m$ as $\mathbb{Z} (K_i)_{\chi_i}$-module if and only if $V_i$ is $FP_m$ as $\mathbb{Z} (G_i)_{\widetilde{\chi}_i}$-module. The proof is completed by the fact that $V_i$ is the trivial $\mathbb{Z} G_i$-module $\mathbb{Z}$. 
\end{proof}

 \begin{proposition} \label{wedgeproducts2} Let $M= \mathbb{Z} X$ be an induced $\mathbb{Z} G$-module and $\chi : G \to \real$ a non-zero  homomorphism. Then the following are equivalent :

1.  $\wedge^i M$ is of type $FP_{m-i}$ as $\z G_{\chi}$-module for all $1 \leq i \leq m$;

2.   $\otimes^i M$ is of type $FP_{m-i}$ as $\z G_{\chi}$-module for all $1 \leq i \leq m$;

3. $\otimes^i M$ is finitely generated as  $\mathbb{Z} G_{\chi}$-module via the diagonal action  and 
 for representatives   $\{ G_1, \ldots , G_{s_1} \}$ of all $G$-orbits of stabilizers of the diagonal action of $G$ on $X^i$,  the restriction $\chi_j$ of $\chi$ on $G_j$ is non-zero and $[\chi_j] \in \Sigma^{m-i} (G_j, \z)$ for all $1 \leq i \leq m$, $1 \leq j \leq s_1$;

4. $\wedge^i M$ is finitely generated as  $\mathbb{Z} G_{\chi}$-module via the diagonal action and 
 for representatives   $\{ \widetilde{G}_1, \ldots , \widetilde{G}_{s_2} \}$ of all $G$-orbits of stabilizers  of the action of $G$ on the $i$-element subsets of $X$, the restriction $\widetilde{\chi}_j$ of $\chi$ on $\widetilde{G}_j$ is non-zero and $[\widetilde{\chi}_j] \in \Sigma^{m-i} (\widetilde{G}_j, \z)$ for all $1 \leq i \leq m$, $1 \leq j \leq s_2$;
 
 5. $G \backslash X^i $ is finite  and for representatives   $\{ G_1, \ldots , G_{s_1} \}$ of all $G$-orbits of stabilizers 
 of the diagonal action of $G$ on $X^i$  the restriction $\chi_j$ of $\chi$ on $G_j$ is non-zero and $[\chi_j] \in \Sigma^{m-i} (G_j, \z)$ for all $1 \leq i \leq m$, $1 \leq j \leq s_1$.
\end{proposition}
 \begin{proof} We say that  the monoid $G_{\chi}$ acts with finitely many orbits on a set $Y$ if there is a finite subset $Y_0$ of $Y$ such that $G_{\chi} Y_0 = Y$. Then  to prove Proposition \ref{wedgeproducts2} it is sufficient to repeat the proof of Proposition
   \ref{wedgeproducts} substituting $G$ with $G_{\chi}$,
  and  apply Lemma \ref{sigmaofactions} and Lemma
   \ref{twisting2} instead of Lemma \ref{twisting}.
 \end{proof}
 
 \begin{theorem} \label{sigmamain}
Let $\Gamma = H \wr_X G$ be a wreath product of type $FP_m$, $X \not= \emptyset$ and $H$ has infinite abelianization.  Let $\chi : \Gamma \to \mathbb{R}$ be a non-zero character such that $\chi(H) = 0$.
Then the following are equivalent :

1. $[\chi]  \in \Sigma^m(\Gamma, \mathbb{Z} )$;

2. $[\chi \mid_G]  \in \Sigma^m(G, \mathbb{Z} )$  and for representatives   $\{ G_1, \ldots , G_{s_1} \}$ of all $G$-orbits of stabilizers 
 of the diagonal action of $G$ on $X^i$  the restriction $\chi_j$ of $\chi$ on $G_j$ is non-zero and $[\chi_j] \in \Sigma^{m-i} (G_j, \z)$ for all $1 \leq i \leq m$, $1 \leq j \leq s_1$. 
\end{theorem}

\begin{proof} The proof is the same as the proof of Theorem
  \ref{mainwreathresult} substituting $G$ with $G_{\chi}$.

The fact that 2. implies 1.  is a monoid version of    Proposition \ref{easydirection}.
  Note that in the proof of Proposition \ref{easydirection}  we
  use that some stabilizers are commensurable. Thus they have the same homological type $FP_k$.
   This step requires justification in its $\Sigma$-version, but this is the content of Lemma \ref{sigmaofactions} (4).

The fact that item 1 implies item 2 follows from Lemma \ref{sigmaretract}
applied to the retracts $G$ and $\Gamma_1 = \mathbb{Z} \wr_X G$ of
$\Gamma$. Indeed we get that $[\chi |_G] \in \Sigma^m(G,
\mathbb{Z})$ and $[\chi |_{\Gamma_1}] \in \Sigma^m(\Gamma_1,
\mathbb{Z})$. By the
monoid version of Theorem \ref{mainthm1} the result holds for $\Gamma_1$, so
condition 5 from Proposition \ref{wedgeproducts2}  holds.

\end{proof}

\end{document}